\documentclass[]{article}
\begin{document}
\newtheorem{proposition}{Proposition}[section]
\newtheorem{definition}{Definition}[section]
\newtheorem{lemma}{Lemma}[section]

\title{\bf New Generalizations of \\ Poisson Algebras }
\author{Keqin Liu\\Department of Mathematics\\The University of British Columbia\\Vancouver, BC\\
Canada, V6T 1Z2}
\date{June 19, 2007}
\maketitle

\begin{abstract} We introduce many new generalizations of Poisson algebras which can be constructed inside the associative algebra of linear transformations over a vector space.
\end{abstract}

\medskip
\section{Introduction}

Poisson algebras play a fundamental role in symplectic geometry. Recently, different generalizations of Poisson algebras have been introduced by several people (\cite{C-K}, 
\cite{C-D}, \cite{K} and \cite{Loday}). A Poisson algebra $P$ is a vector space equipped with a square bracket $[\, ,\,]$ and a dot product $\cdot$ such that the following three conditions hold: 
\begin{description}
\item[Condition 1] $P$ is a Lie algebra with respect to the square bracket $[\, ,\,]$;
\item[Condition 2] $P$ is a commutative associative algebra with respect to the dot product $\cdot$\,\,;
\item[Condition 3] The square bracket $[\, ,\,]$ and the dot product $\cdot$ satisfy the Leibniz rule: $[ x, y\cdot z]=[x, y]\cdot z+y\cdot [x, z]$ for $x$, $y$ and $z\in P$.
\end{description}
In this paper, we introduce many new generalizations of Poisson algebras. All of the generalizations of Poisson algebras can be divided into two groups. One group consists of the tailless Poisson algebras which do not use a derivative to generalize the Leibniz rule. The other group consists of the tailed Poisson algebras which use a derivative to generalize the Leibniz rule. There are a few kinds of tailless Poisson algebras, but there are many kinds of tailed Poisson algebras. The details of constructing these generalizations of Poisson algebras will be given in \cite{Liu}. The important facts about these generalizations are that each of these generalizations of Poisson algebras can be constructed inside the associative algebra of linear transformations over a vector space, and some of  these generalizations of Poisson algebras can be constructed inside the associative algebra of linear transformations over a vector space by using many different ways.

\section{Tailless Poisson Algebras}

After dropping off the commutative property in the Condition 2, we get the following

\begin{definition}\label{def2.2} A {\bf square-circle Poisson algebra} $P$ is a vector space equipped with a {\bf square bracket} $[\, , \,]$ and a {\bf circle product} $\circ$ which have three properties:
\begin{description}
\item[(i)] $P$ is a Lie algebra with respect to the square bracket $[\, ,\,]$;
\item[(ii)] $P$ is an associative algebra with respect to the circle product $\circ$;
\item[(iii)] The square bracket $[\, , \,]$ and the circle product $\circ$ satisfy the 
{\bf Leibniz rule}:
\begin{equation}
[x,  y\circ z]=[x, y]\circ z +y\circ [x, z]\qquad\mbox{for $x,y,z\in P$.}
\end{equation}
\end{description}
\end{definition}

A square-circle Poisson algebra is called a non-commutative Poisson algebra in \cite{K}.
Recall from \cite{Loday} that a {\bf right Leibniz algebra} $L$ is a vector space equipped with a {\bf angle bracket} $\langle \,, \,\rangle$ satisfying the {\bf right Leibniz identity}; that is, 
\begin{equation}
\langle x, \langle y, z\rangle  \rangle =\langle\langle x,  y\rangle , z\rangle 
 - \langle\langle x,  z\rangle , y\rangle \qquad\mbox{for $x,y,z\in L$.}
\end{equation}

\begin{definition} Let $P$ be a vector space equipped with a {\bf angle bracket}
$\langle \, , \,\rangle$ and a {\bf circle product} $\circ$ such that  $P$ is a right Leibniz algebra with respect to the angle bracket $\langle \, , \,\rangle$ and $P$ is an associative algebra with respect to the circle product $\circ$.
\begin{description}
\item[(i)] $P$ is called a  {\bf left angle-circle Poisson algebra} if the angle bracket 
$\langle \, , \,\rangle$ and the circle product $\circ$ satisfy the {\bf left Leibniz rule}:
\begin{equation}
\langle x,  y\circ z\rangle=\langle x, y\rangle\circ z +y\circ \langle x, z\rangle
\qquad\mbox{for $x,y,z\in P$.}
\end{equation}
\item[(ii)] $P$ is called a  {\bf right angle-circle Poisson algebra} if the angle bracket 
$\langle \, , \,\rangle$ and the circle product $\circ$ satisfy the {\bf right Leibniz rule}:
\begin{equation}
\langle x\circ y , z\rangle =x\circ \langle y, \,z\rangle + \langle x, \,z\rangle \circ y
\qquad\mbox{for $x,y,z\in P$.}
\end{equation}
\end{description}
\end{definition}

A right angle-circle Poisson algebra is called a non-commutative Leibniz-Poisson algebra in \cite{C-D}.

\medskip
The next proposition gives a summary of our ways of constructing  tailless Poisson algebras.

\begin{proposition}\label{pr1} Let $End (V)$ be the associative algebra of linear transformations over a 
vector space $V$. If $V$ has a proper nonzero subspace, then there exists a subspace $\mathcal{M}$ of $End(V)$ such that the following hold:
\begin{description}
\item[(i)] There are $46$ ways of making $\mathcal{M}$ into a square-circle Poisson algebra;
\item[(ii)] There are $6$ ways of making $\mathcal{M}$ into a left angle-circle Poisson algebra;
\item[(iii)] There are $51$ ways of making $\mathcal{M}$ into a right angle-circle Poisson algebra.
\end{description}
\end{proposition}

A remark about left angle-circle Poisson algebras is that each of the $6$ ways of making $\mathcal{M}$ into a left angle-circle Poisson algebra in Proposition~\ref{pr1} simultaneously makes $\mathcal{M}$ into a right angle-circle Poisson algebra. 

\section{Tailed Square-Circle Poisson Algebras}

We begin this section with the following

\begin{definition}\label{def1.2} Let $A$ be an associative algebra with respect to an associative product $\circ$, and let $x$, $y$ and $z$ are elements of $A$.
\begin{description}
\item[(i)] If $A$ is a Lie algebra with respect to a square bracket $[\, , \,]$, then \begin{equation}\label{eq1.7}
J[x, y, \circ, z]:=[x, \, y \circ z] -[x, \, y] \circ z-y\circ [x, \, z] 
\end{equation}
is called the {\bf square-circle Jacobian}.
\item[(ii)] If $A$ is a right Leibniz algebra with respect to a angle bracket 
$\langle \, , \,\rangle $, then 
\begin{equation}\label{eq1.8}
J_{\ell}\langle x, y, \circ , z\rangle := \langle x, \,y\circ z\rangle 
-\langle x, \,y\rangle \circ z-y\circ \langle x, \,z\rangle  
\end{equation}
is called the {\bf left angle-circle Jacobian}, and
\begin{equation}\label{eq1.9}
J_{r}\langle x, \circ, y, z\rangle := \langle x\circ y , z\rangle 
-x\circ \langle y, \,z\rangle - \langle x, \,z\rangle \circ y 
\end{equation}
is called the {\bf right angle-circle Jacobian}.
\end{description}
\end{definition}

We now introduce $8$ different tailed square-circle Poisson algebras in the following

\begin{definition}\label{def1.3} Let $P$ be a vector space equipped with a {\bf square bracket} 
$[ \, , \, ]$ and a {\bf circle product} $\circ$ such that  $P$ is a Lie algebra with respect to the square bracket $[ \, , \, ]$ and $P$ is an associative algebra with respect to the circle product $\circ$. Let $D$ be a derivation with respect to both the square bracket $[\, , \,]$ and the circle product $\circ$. Let $x$, $y$ and $z$ be arbitrary three elements of $P$.
\begin{description}
\item[(i)] $P$ is called a {\bf tailed\,$^{1-st}$ square-circle Poisson algebra} if
\begin{equation}\label{eq1.10}
J[x, y, \circ, z]=x\circ y\circ D(z)-y\circ x\circ D(z) .
\end{equation}
\item[(ii)] $P$ is called a {\bf tailed\,$^{2-nd}$ square-circle Poisson algebra} if
\begin{equation}\label{eq1.19}
J[x, y, \circ, z]=y\circ D(x)\circ z .
\end{equation}
\item[(iii)] $P$ is called a {\bf tailed\,$^{3-rd}$ square-circle Poisson algebra} if there exists $\alpha\in \mathbf{k}$ such that
\begin{equation}\label{eq1.34}
J[x, y, \circ, z]=D(y)\circ z\circ x-D(y)\circ x\circ z+ y\circ (\alpha D)(x)\circ z .
\end{equation}
\item[(iv)] $P$ is called a {\bf tailed\,$^{4-th}$ square-circle Poisson algebra} if
\begin{equation}\label{eq1.43}
J[x, y, \circ, z]=D(y)\circ z\circ x-D(y)\circ x\circ z
\end{equation}
and
\begin{equation}\label{eq1.44}
x\circ y\circ D(z)=x\circ D(y)\circ z=0 .
\end{equation}
\item[(v)] $P$ is called a {\bf tailed\,$^{5-th}$ square-circle Poisson algebra} if there exists a nonzero scalar $\alpha\in \mathbf{k}$ such that
\begin{equation}\label{eq1.76}
J[x, y, \circ, z]=y\circ (\alpha D)(x)\circ z+D(y)\circ z\circ x-D(y)\circ x\circ z .
\end{equation}
\item[(vi)] $P$ is called a {\bf tailed\,$^{6-th}$ square-circle Poisson algebra} if there exists a scalar $\alpha\in \mathbf{k}$ such that
\begin{eqnarray}\label{eq1.81}
J[x, y, \circ, z]&=&(\alpha D)(y)\circ z\circ x-(\alpha D)(y)\circ x\circ z +\nonumber\\
&&\quad +x\circ y\circ D(z)-y\circ x\circ D(z) .
\end{eqnarray}
\item[(vii)] $P$ is called a {\bf tailed\,$^{7-th}$ square-circle Poisson algebra} if there exists a nonzero scalar $\alpha\in \mathbf{k}$ such that
\begin{equation}\label{eq1.86}
J[x, y, \circ, z]=y\circ (\alpha D)(x)\circ z +x\circ y\circ D(z)-y\circ x\circ D(z) .
\end{equation}
\item[(viii)] $P$ is called a {\bf tailed\,$^{8-th}$ square-circle Poisson algebra} if 
\begin{equation}\label{eq1.91}
J[x, y, \circ, z]=x\circ y\circ D(z)-y\circ x\circ D(z)
\end{equation}
and
\begin{equation}\label{eq1.92}
D(x)\circ y\circ z=x\circ D(y)\circ z=0 .
\end{equation}
\end{description}
\end{definition}

\medskip
Note that if the circle product $\circ$ is commutative, then the tail $y\circ D(x)\circ z$ of a tailed$^{2-nd}$ square-circle Poisson algebra becomes the tail $D(x)\circ y\circ z$ of a generalized Poisson superalgebra introduced in \cite{C-K}.

\medskip
The next proposition gives a summary of our ways of constructing   
tailed square-circle Poisson algebras.

\begin{proposition}\label{pr2} Let $End (V)$ be the associative algebra of linear transformations over a vector space $V$. If $V$ has a proper nonzero subspace, then there exists a subspace $\mathcal{M}$ of $End(V)$ such that the following hold:
\begin{description}
\item[(i)] There are $2$ ways of making $\mathcal{M}$ into a tailed\,$^{1-st}$ square-circle Poisson algebra;
\item[(ii)] There are $4$ ways of making $\mathcal{M}$ into a tailed\,$^{2-nd}$ square-circle Poisson algebra;
\item[(iii)] There are $2$ ways of making $\mathcal{M}$ into a tailed\,$^{3-rd}$ square-circle Poisson algebra;
\item[(iv)] There are $7$ ways of making $\mathcal{M}$ into a tailed\,$^{4-th}$ square-circle Poisson algebra;
\item[(v)] There are $1$ way of making $\mathcal{M}$ into a tailed\,$^{i-th}$ square-circle Poisson algebra for $i=5$, $6$ and $7$;
\item[(vi)] There are $6$ ways of making $\mathcal{M}$ into a tailed\,$^{8-th}$ square-circle Poisson algebra.
\end{description}
\end{proposition}

\section{Tailed Left Angle-Circle Poisson Algebras}

In this section, we introduce $6$ different tailed left angle-circle Poisson algebras. Their definitions are as follows.

\begin{definition} Let $P$ be a vector space equipped with a {\bf angle bracket} 
$\langle \, , \, \rangle$ and a {\bf circle product} $\circ$ such that  $P$ is a right Leibniz algebra with respect to the angle bracket $\langle \, , \, \rangle$ and $P$ is an associative algebra with respect to the circle product $\circ$. Let $D$ be a derivation with respect to both the angle bracket $\langle \, , \, \rangle$ and the circle product $\circ$. Let $x$, $y$ and $z$ be arbitrary three elements of $P$.
\begin{description}
\item[(i)] $P$ is called a {\bf tailed\,$^{1-st}$ left angle-circle Poisson algebra} if
\begin{equation}
J_{\ell}\langle x, y, \circ , z\rangle=x\circ y\circ D(z)-y\circ x\circ D(z) .
\end{equation}
\item[(ii)] $P$ is called a {\bf tailed\,$^{2-nd}$ left angle-circle Poisson algebra} if
\begin{equation}
J_{\ell}\langle x, y, \circ , z\rangle=x\circ y\circ D(z)-y\circ x\circ D(z)
\end{equation}
and 
\begin{equation}
x\circ D(y)\circ z=D(x)\circ y\circ z=0.
\end{equation}
\item[(iii)] $P$ is called a {\bf tailed\,$^{3-rd}$ left angle-circle Poisson algebra} if
\begin{equation}
J_{\ell}\langle x, y, \circ , z\rangle=D(y)\circ x\circ z-D(y)\circ z\circ x.
\end{equation}
\item[(iv)] $P$ is called a {\bf tailed\,$^{4-th}$ left angle-circle Poisson algebra} if
\begin{equation}
J_{\ell}\langle x, y, \circ , z\rangle=D(y)\circ x\circ z-D(y)\circ z\circ x
\end{equation}
and
\begin{equation}
x\circ y\circ D(z)=x\circ D(y)\circ z .
\end{equation}
\item[(v)] $P$ is called a {\bf tailed\,$^{5-th}$ left angle-circle Poisson algebra} if
\begin{eqnarray}
J_{\ell}\langle x, y, \circ , z\rangle&=&
x\circ y\circ  D(z)-y\circ x\circ D(z)+\nonumber\\
&& -D(y)\circ x\circ z+D(y)\circ z\circ x
\end{eqnarray}
\item[(vi)] $P$ is called a {\bf tailed\,$^{6-th}$ left angle-circle Poisson algebra} if
\begin{eqnarray}
J_{\ell}\langle x, y, \circ , z\rangle&=&x\circ y\circ D(z)-
y\circ x\circ D(z)
\end{eqnarray}
and
\begin{eqnarray}
&&x\circ D(y)\circ z=D(x)\circ y\circ z .
\end{eqnarray}
\end{description}
\end{definition}

\medskip
The next proposition gives a summary of our ways of constructing  tailed left angle-circle Poisson algebras.

\begin{proposition} Let $End (V)$ be the associative algebra of linear transformations over a vector space $V$. If $V$ has a proper nonzero subspace, then there exists a subspace $\mathcal{M}$ of $End(V)$ such that the following hold:
\begin{description}
\item[(i)] There are $3$ ways of making $\mathcal{M}$ into a tailed\,$^{1-st}$ left angle-circle Poisson algebra;
\item[(ii)] There are $5$ ways of making $\mathcal{M}$ into a tailed\,$^{2-nd}$ left angle-circle Poisson algebra;
\item[(iii)] There are $2$ ways of making $\mathcal{M}$ into a tailed\,$^{3-rd}$ left angle-circle Poisson algebra;
\item[(iv)] There are $6$ ways of making $\mathcal{M}$ into a tailed\,$^{4-th}$ left angle-circle Poisson algebra;
\item[(v)] There are $1$ way of making $\mathcal{M}$ into a tailed\,$^{i-th}$ left angle-circle Poisson algebra for $i=5$ and $6$.
\end{description}
\end{proposition}

\section{Tailed Right Angle-Circle Poisson Algebras}

We now introduce $4$ different tailed right angle-circle Poisson algebras. 

\begin{definition} Let $P$ be a vector space equipped with a {\bf angle bracket} 
$\langle \, , \, \rangle$ and a {\bf circle product} $\circ$ such that  $P$ is a right Leibniz algebra with respect to the angle bracket $\langle \, , \, \rangle$ and $P$ is an associative algebra with respect to the circle product $\circ$. Let $D$ be a derivation with respect to both the angle bracket $\langle \, , \, \rangle$ and the circle product $\circ$. Let $x$, $y$ and $z$ be arbitrary three elements of $P$.
\begin{description}
\item[(i)] $P$ is called a {\bf tailed\,$^{1-st}$ right angle-circle Poisson algebra} if
\begin{equation}\label{eq3.1}
J_{r}\langle x, \circ, y, z\rangle=x\circ z\circ D(y)-z\circ x\circ D(y) .
\end{equation}
\item[(ii)] $P$ is called a {\bf tailed\,$^{2-nd}$ right angle-circle Poisson algebra} if
\begin{equation}
J_{r}\langle x, \circ, y, z\rangle=x\circ z\circ D(y)-z\circ x\circ D(y)
\end{equation}
and
\begin{equation}
x\circ D(y)\circ z=D(x)\circ y\circ z .
\end{equation}
\item[(iii)] $P$ is called a {\bf tailed\,$^{3-rd}$ right angle-circle Poisson algebra} if
\begin{equation}
J_{r}\langle x, \circ, y, z\rangle=D(x)\circ y\circ z-D(x)\circ z\circ y .
\end{equation}
\item[(iv)] $P$ is called a {\bf tailed\,$^{4-th}$ right angle-circle Poisson algebra} if
\begin{equation}
J_{r}\langle x, \circ, y, z\rangle=D(x)\circ y\circ z-D(x)\circ z\circ y
\end{equation}
and
\begin{equation}
x\circ y\circ D(z)=x\circ D(y)\circ z=0 .
\end{equation}
\end{description}
\end{definition}

\medskip
The next proposition gives a summary of our ways of constructing  tailed right angle-circle Poisson algebras.

\begin{proposition} Let $End (V)$ be the associative algebra of linear transformations over a vector space $V$. If $V$ has a proper nonzero subspace, then there exists a subspace $\mathcal{M}$ of $End(V)$ such that the following hold:
\begin{description}
\item[(i)] There are $2$ ways of making $\mathcal{M}$ into a tailed\,$^{1-st}$ right angle-circle Poisson algebra;
\item[(ii)] There are $3$ ways of making $\mathcal{M}$ into a tailed\,$^{2-nd}$ right angle-circle Poisson algebra;
\item[(iii)] There are $2$ ways of making $\mathcal{M}$ into a tailed\,$^{3-rd}$ right angle-circle Poisson algebra;
\item[(iv)] There are $3$ ways of making $\mathcal{M}$ into a tailed\,$^{4-th}$ right angle-circle Poisson algebra.
\end{description}
\end{proposition}

\end{document}